\documentstyle[newamssymbols]{birk-ot}
\Year{2001}
\newcommand{\bH}{{H}}

\newcommand{\cB}{{\mathcal B}}
\newcommand{\cD}{{\mathcal D}}

\newcommand{\cH}{{\mathcal H}}

\newcommand{\cO}{{\mathcal O}}

\newcommand{\ri}{{\rm i}}

\newcommand{\sP}{{\sf P}}

\newcommand{\tA}{\widetilde{A}}
\newcommand{\ta}{\widetilde{a}}
\newcommand{\tB}{\widetilde{B}}

\newcommand{\tD}{\widetilde{D}}
\newcommand{\trl}{\triangleleft}
\newcommand{\trr}{\triangleright}

\newcommand{\vkappa}{\varkappa}
\newcommand{\D}{\displaystyle}
\newcommand{\dom}{\mathop{\rm dom}}

\newcommand{\Inf}{\mathop{\rm inf}}
\newcommand{\Int}{\displaystyle\int\limits}
\newcommand{\Sup}{\mathop{\rm sup}}

\newcommand{\dist}{\mathop{\rm dist}}
\newcommand{\Img}{\mathop{\rm Im}}
\newcommand{\Real}{\mathop{\rm Re}}

\newcommand{\lal}{{\langle}}
\newcommand{\ral}{{\rangle}}
\newcommand{\Var}{{\rm Var}}
\newcommand{\acs}{{\sigma_{\rm ac}}}

\newcommand{\spec}{\mathop{\sigma}}
\newcommand{\ran}{\mathop{\rm ran}}
\newcommand{\res}{\mathop{\varrho}}

\newtheorem{hypothesis}[theorem]{Hypothesis}
\begin{document}
\Title{Factorization Theorem for the Transfer Function
Associated with an Unbounded Non-Self-Adjoint
${\bf 2}\times{\bf 2}$ Operator Matrix}
\Shorttitle{Factorization Theorem}
\By{{\sc V. Hardt, R. Mennicken, A. K. Motovilov}}
\Names{Hardt, Mennicken, Motovilov}
\Email{motovilv@thsun1.jinr.ru}
\maketitle
%
\begin{abstract}
We construct operators which factorize the
transfer function associated with a non-self-adjoint
\mbox{$2\times2$} operator matrix whose diagonal entries can have
overlapping spectra and whose off-diagonal entries are unbounded
operators.
\end{abstract}

\newsection{Introduction}\label{Intro}
In the present work we extend some of the results of Refs.
\cite{HMM1,HMM2} obtained for a class of unbounded
self-adjoint block operator matrices to the case where these
matrices are unbounded non-self-adjoint operators.

Throughout the paper we consider a $2\times 2$ block operator matrix
having the form
\begin{equation}
\label{twochannel}
\bH_0=
\left(\begin{array}{ccc}
A    & &  B \\
D    & &  C     \end{array}\right)
\end{equation}
and acting in the orthogonal sum  $\cH=\cH_A\oplus\cH_C$ of separable
Hilbert spaces $\cH_A$ and $\cH_C$.  The entry $C$
is a possibly unbounded self-adjoint
operator in $\cH_C$ on the domain  $\dom(C)$. This operator is assumed
to be semibounded from below, that is,
\begin{equation}
\label{Csemib}
C\geq\lambda_C \quad\mbox{for some}\quad \lambda_C\in\R.
\end{equation}
In the following we assume without loss of generality that the
lower bound $\lambda_C$ for the entry $C$ in the assumption
(\ref{Csemib}) is positive,
$$
\lambda_C>0
$$
(otherwise one could simply shift the origin of the spectral
parameter axis).  The entry $A$ is supposed to be a bounded and
not necessarily self-adjoint operator in $\cH_A$, that is,
$A\in\cB(\cH_A)$.

Regarding the couplings $B$ and $D$ we assume the following
hypotheses:
\begin{enumerate}
\item[(i)] {\it $D$ is a densely defined closable operator
from $\dom(D)\subset\cH_A$ to $\cH_C$ such that
the product $C^{-1/2}D$ admits extension from $\dom(D)$ to
the whole space $\cH_A$ as a bounded operator in
$\cB(\cH_A,\cH_C)$,}
\item[(ii)] {\it $B$ is a densely defined closable operator
from $\dom(B)\subset\cH_C$ to $\cH_A$} such that
\end{enumerate}
\begin{equation}
\label{IniCond}
\dom(B)\supset\dom(C^{1/2}).
\end{equation}
Since $\dom(C^{1/2})\supset\dom(C)$, these assumptions imply
that the matrix $\bH_0$ is a densely defined closable operator
on the domain $\dom(\bH_0)=\dom(D)\oplus\dom(C)$. We denote
its closure by $H$, that is, $\bH=\overline{\bH}_0$.

From the assumption (i) one infers that
$$
\tD:=\overline{C^{-1/2}D}\in\cB(\cH_A,\cH_C).
$$
Also one concludes from (i) that
$$
\ran\left(\left.\tD \right|_{\dom(D)}\right)\subset\dom(C^{1/2})
$$
and
$$
 D=\left.{C^{1/2} \tD}\right|_{\dom(D)}.
$$
At the same time the assumption (\ref{IniCond}) yields that the
product $\tB:=BC^{-1/2}$ is a bounded linear operator between
$\cH_C$ and $\cH_A$, i.\,e. $\tB\in\cB(\cH_C,\cH_A)$, and
$\left.B\right|_{\dom(C^{1/2})}=\tB C^{1/2}$.

In addition, the hypotheses (i) and (ii) imply that
for $z$ in the resolvent set $\res(C)$ of $C$
the operator-valued function $A-z-B(C-z)^{-1}D$ is
densely defined on $\dom(D)$ and admits
a bounded extension onto the whole space
$\cH_A$. We denote this extension by $M(z)$,
\begin{equation}
\label{transferfunction}
M(z):=A-z-\overline{B(C-z)^{-1}D},
\end{equation}
and call it the transfer function associated with the operator $\bH$.
By definition
\begin{eqnarray*}
\left.{M(z)}\right|_{\dom(D)} &= &
A-z - B(C-z)^{-1}D\\
&=&A-z - \tB C^{1/2}(C-z)^{-1}C^{1/2}\tD\\
&=& \tA-z+V(z),
\end{eqnarray*}
where
\begin{equation}
\label{tA-V}
\tA:=A-\tB\tD \quad \mbox{ and } \quad
V(z):=z\,\tB(z-C)^{-1}\tD.
\end{equation}
Meanwhile $\tA$ is a bounded operator on $\cH_A$ while $V(z)$,
$z\in\res(C)$, is an operator-valued function with values in
$\cB(\cH_A)$. Then it follows that
\begin{equation}
\label{Mz}
M(z)=\tA-z+V(z).
\end{equation}

Before going into more details we mention the recent
paper~\cite{LMMT} studying by means of the new concept of the
quadratic numerical range, introduced in~\cite{LT}, in
particular the factorization properties of the transfer
functions associated with bounded non-self-adjoint $2\times2$
block operator matrices. As a matter of fact, the operators
factorizing the transfer functions for bounded non-self-adjoint
operator matrices have been also constructed in \cite{AdLT} and
for unbounded non-self-adjoint operator matrices with bounded
entries $B$ and $D$ in \cite{LT}.  The results of
\cite{AdLT,LMMT,LT}  are obtained under the assumption that the
spectra of $A$ and $C$ do not overlap.

We are studying the transfer function (\ref{Mz}), associated
with the unbounded non-self-adjoint operator
matrix (\ref{twochannel}) satisfying the assumptions (i) and
(ii), in the case where the numerical range $\nu(\tA)$ of the
operator $\tA$ (and in particular its spectrum $\spec(A)$) may
have a non-empty intersection with the spectrum $\spec(C)$ of
$C$. Notice that, since the resolvent of the operator $\bH$ can
be expressed explicitly in terms of $\bigl[M(z)\bigr]^{-1}$
(see, e.\,g., \cite{AdLMSr,MenShk,MennMotMathNachr}), in
studying the spectral properties of the transfer function $M$
one studies at the same time the spectral properties of the
operator matrix $\bH$.

Obviously, the transfer function (\ref{Mz}), considered for
$z\in\res(C)$, represents a holomorphic operator-valued
function.  (We refer to the standard definition of holomorphy of
an operator-valued function with respect to the operator norm
topology, see, e.\,g., \cite{AdLMSr}.) We study the transfer
function $M(z)$ under the assumption that it admits analytic
continuation through the absolutely continuous spectrum of the
entry $C$ at least in a neighborhood of the numerical range
$\nu(\tA)$ of the operator $\tA$.  In order to present our main
ideas in a more transparent form, we suppose, for the sake of
simplicity, that the entry $C$ only has the absolutely
continuous spectrum, that is,
$$
  \spec(C)=\acs(C).
$$
Section\,2 includes a description of the conditions making the
analytic continuation of $M(z)$ through the set $\acs(C)$ to be
possible.  In Section\,3 we introduce the nonlinear
transformation equations (\ref{MainEq}) making a rigorous sense
to the formal operator equation $M(Z)=0$, $Z\in\cB(\cH_A)$. We
explicitly show that eigenvalues and accompanying eigenvectors
of a solution $Z$ to the equation (\ref{MainEq}) are eigenvalues
and eigenvectors of the analytically continued transfer function
$M$. Further, we prove the solvability of (\ref{MainEq}) under
smallness conditions concerning the operator $\tB$ and $\tD$, see
(\ref{Best}).  In Section\,4 a factorization theorem (Theorem
\rref{factorization}) is proven for the analytically continued
transfer function. This theorem implies in particular that there
exist certain domains in $\C$ lying partly on the unphysical
sheet(s) where the spectrum of the analytically continued
transfer function is represented by the spectrum of the
corresponding solutions to the transformation equations
(\ref{MainEq}).  Finally, in Section\,5 we present an example.

\newsection{Analytic continuation}
\label{Transfer_function}
Throughout this paper we assume that the spectrum $\spec(C)$ of
the entry $C$ is absolutely continuous and it fills the semiaxis
$[\lambda_C,+\infty)$.  By $\nu(\tA)$ we denote the numerical
range of the (bounded) operator $\tA$,
$$
\nu(\tA)=\overline{\left\{\lal\tA x,x\ral:\, x\in\cH_A,\,\|x\|=1\right\}}.
$$
We suppose that the set $\nu(\tA)$ lyes
inside a neighborhood
$$
\cO_\eta([\alpha_1,\alpha_2]):=\{z\in\C:\,
\dist(z,[\alpha_1,\alpha_2])\leq\eta\}, \quad \eta>0,
$$
of a finite real interval $[\alpha_1,\alpha_2]\subset\R$,
$\alpha_1<\alpha_2$, and $\alpha_1-\eta>\lambda_C$.  Notice that
the numerical range $\nu(\tA)$ is a convex set containing the
spectrum $\spec(\tA)$ (see, e.\,g.,~\cite{GK}, \S\,V.6) and,
moreover,
\begin{equation}
\label{Arest}
\|(\tA-zI_A)^{-1}\|\leq\frac{1}{\dist\big(z,\nu(\tA)\big)},\quad
z\in\C\setminus\nu(\tA),
\end{equation}
where $I_A$ stands for the identity operator in $\cH_A$
(see Lemma~V.6.1 in~\cite{GK}).

Let $\{E_C(\mu)\}_{\mu\in\R}$ be the spectral family for the
entry $C$, $C=\int_{\spec(C)}\mu\,dE_C(\mu)$.  Then the
function $V(z)$ can be written
$$
    V(z) = \int_{\lambda_C}^\infty dK(\mu)\frac{z}{z-\mu}
$$
with
$$
  K(\mu) := \tB\widehat{E}_C(\mu)\tD
$$
We assume that the function $K(\mu)$ is continuously
differentiable in $\mu\in(\lambda_C,+\infty)$ in the operator
norm topology and, moreover, that it admits analytic
continuation from the interval $(\lambda_C,\beta)$,
$\beta>\alpha_2+\eta$, to a simply connected domain
$\cD\subset\C$, $\cD\supset\cO_\eta([\alpha_1,\alpha_2])$.  For
the continuation we keep the same notation $K(\mu)$ and by
$K'(\mu)$ denote the derivative of $K$.  We
suppose that at the points $\lambda=\lambda_C$ and
$\lambda=\beta$ the function $K'(\mu)$ satisfies the condition
$$
\|K'(\mu)\|\leq c|\mu-\lambda|^\gamma, \quad\mu\in \cD,
$$
with some $c>0$ and $\gamma\in (-1,0]$. We also assume that the
operators $\tB$ and $\tD$
are such that
\begin{equation}
\label{NBnorm}
\displaystyle\int_\beta^{\infty}|d\mu|\,(1+|\mu|)^{-1}\|K'(\mu)\|<\infty.
\end{equation}

In the following we use the notation
$$\cD^{(+1)}:=\cD\cap\C^+
\quad
\mbox{and}\quad \cD^{(-1)}:=\cD\cap\C^-.
$$
\begin{lemma}
\label{M1-Continuation}
Let $\widetilde{\Gamma}_l$ $(l=\pm 1)$ be a rectifiable Jordan
curve in $\cD^{(l)}$ resulting from continuous deformation of the
interval $(\lambda_C,\beta)$, the end points of this
interval being fixed, and let
${\Gamma}_{l}=\widetilde{\Gamma}_{l}\cup[\beta,+\infty)$.
Then the analytic continuation of the transfer function $M(z)$,
$z\in\C\setminus[\lambda_C,+\infty)$, through the spectral interval
$(\lambda_C,\beta)$ into the subdomain $\cD(\Gamma_l)\subset
\cD^{(l)}$ $(l=\pm1)$ bounded by $(\lambda_C,\beta)$ and
$\widetilde{\Gamma}_l$ is given by
\begin{equation}\label{Mcmpl}
 M_{\Gamma_l}(z) := \tA-z+V_{\Gamma_l}(z),
\end{equation}
where
\begin{equation}
\label{MGamma}
V_{\Gamma_l}(z) :=
\int_{\Gamma_l} d\mu\,K'(\mu)\,\frac{z}{z-\mu}.
\end{equation}
For $z\in\cD(\Gamma_l)$ the function
$M_{\Gamma_l}(z)$ may be written as
\begin{equation}\label{M1Gresidue}
  M_{\Gamma_l}(z) = M(z)+2\pi\ri\,l z K'(z).
\end{equation}
\end{lemma}
\begin{proof}
The function~(\ref{MGamma}) is well defined for $z\not\in\Gamma_l$
since (\ref{NBnorm}) holds and for
any $z\in\C\setminus \Gamma_l$ there exist a number $c(z)>0$ such that
the estimate $\bigl |(z-\mu)^{-1}\bigr |<c(z)\,\bigl
(1+|\mu|\bigr )^{-1}$ $(\mu\in \Gamma _l\bigr )$ is valid.  Then
the proof of this lemma is reduced to the observation that the
function $M_{\Gamma_l}(z)$ is holomorphic for $z\in\C\setminus
\Gamma_l$ and coincides with $M(z)$ for $z\in\C\setminus
\overline{\cD(\Gamma_l)}$.  The equation~(\ref{M1Gresidue}) is
obtained from~(\ref{MGamma}) using the Residue Theorem.
\end{proof}

\begin{remark}
{}From formula (\ref{M1Gresidue}) one concludes that
the transfer function $M(z)$ has a Riemann surface larger
than a single sheet of the spectral parameter plane.
The sheet of the complex plane where the transfer
function $M(z)$ together with the resolvent $(\bH -z)^{-1}$ is
initially considered  is said to be the {physical sheet}.
Hence, the formula (\ref{M1Gresidue}) implies that the domains
$\cD^{(+1)}$ and $\cD^{(-1)}$ should be placed on additional sheet(s)
of the complex plane different from the physical
sheet. Recall that these additional sheets are usually called
unphysical sheets (see, e.\,g., \cite{ReedSimonIII}).
\end{remark}

\begin{remark}
\label{r2.2}
For $z\in\C\setminus\Gamma_l$, the  formula {\rm (\ref{MGamma})}
defines values of the function $V_{\Gamma_l}(\cdot)$ in
$\cB(\cH_A)$.  The inverse transfer function $\bigl[M(z)\bigr
]^{-1}$ coincides with the compressed resolvent
$P_{\cH_A}\left.(\bH-z)^{-1}\right|_{\cH_A}$ where $P_{\cH_A}$
stands for the orthogonal projection on $\cH_A$.  Thus,
$\bigl[M(z)\bigr ]^{-1}$ is holomorphic in
$\C\setminus\spec(\bH)$. Since $M_{\Gamma_l}(z)$ coincides with
$M(z)$ for all $z\in\C\setminus \overline{\cD(\Gamma_l)}$, one
concludes that $[M_{\Gamma_l}(z)]^{-1}$ exists as a bounded
operator and is holomorphic in $z$ at least for
$z\in\big(\C\setminus\spec(\bH)\cup\overline{\cD(\Gamma_l)}$.
\end{remark}

\newsection{The transformation equations}
\label{SmainEq}
Let $\Gamma_l\subset\cD^{(l)}$ $(l=\pm1)$ be a contour described
in hypothesis of Lemma \ref{M1-Continuation}.  Assume that
$Z\in\cB(\cH_A)$ is a bounded operator such that its spectrum
$\spec(Z)$ is separated from $\Gamma_l$,
$$
\dist\big(\spec(Z),\Gamma_l\big)>0.
$$
Following to
\cite{MennMotMathNachr,MotSPbWorkshop,MotRem} (cf.~\cite{AdLT}),
for such $Z$ we introduce the ``right'',
$V^\trr_{\Gamma_l}(Z)$, and ``left'', $V^\trl_{\Gamma_l}(Z)$, transformations
respectively
\begin{equation}
\label{V1Y}
V^\trr_{\Gamma_l}(Z) := \int_{\Gamma_l}
d\mu\,K'(\mu)\,Z(Z-\mu)^{-1}.
\end{equation}
and
\begin{equation}
\label{V1Yl}
V^\trl_{\Gamma_l}(Z) := \int_{\Gamma_l}
Z(Z-\mu)^{-1}K'(\mu)\,d\mu.
\end{equation}
Obviously, for both ``right'' and ``left''
symbols $\varkappa=\trr$ and $\varkappa=\trl$
\begin{eqnarray}
\nonumber
\|V^\varkappa_{\Gamma_l}(Z)\| &\leq& \|Z\| \Sup_{\mu\in\Gamma_l}
\left[(1+|\mu|)\|(Z-\mu)^{-1}\|\right]\\
\label{V1Yest}
 && \times
\displaystyle\int_{\Gamma_l}|d\mu|\,(1+|\mu|)^{-1}\|K'(\mu)\|<\infty.
\end{eqnarray}

In what follows we consider the operator transformation equations on
$\cB(\cH_A)$ (cf. \cite{AdLT,MennMotMathNachr,MotSPbWorkshop,MotRem})
\begin{equation}
\label{MainEq}
Z = \tA+V^\vkappa_{\Gamma_l}(Z),\quad
  \vkappa=\trr\quad\mbox{or}\quad\vkappa=\trl\,.
\end{equation}
In particular, for $\vkappa=\trr$ the equation
(\ref{MainEq}) possesses the following
characteristic property: If an operator $Z_\trr$ is a
solution of~(\ref{MainEq}) and $u\in\cH_A$ is an eigenvector of $Z_\trr$,
i.\,e., $Z_\trr u=zu$ for some $z\in\spec(Z_\trr)$, then
\begin{eqnarray*}
zu&=&\tA u+V^\trr_\Gamma(Z_\trr)u =
\tA u+\int_{\Gamma_l}d\mu\,K'(\mu )Z_\trr(Z_\trr-\mu)^{-1}u\\
&=&\tA u+\int_{\Gamma_l}d\mu\,K'(\mu )\frac z{z-\mu}u
 = \tA u+V_{\Gamma_l}(z)u.
\end{eqnarray*}
Hence, any eigenvalue $z$ of such an operator $Z_\trr$ is
automatically an eigenvalue for the analytically continued
transfer function $M_{\Gamma_l}(\cdot)$ and $u$ is a corresponding
eigenvector. One can easily see that a similar relation holds
between the operator $Z_\trl^*$, adjoint of a solution $Z_\trl$
to the transformation equation~(\ref{MainEq}) for
$\vkappa=\trl$, and the adjoint transfer function
$[M_{\Gamma_l}(\cdot)]^*$.  This means that, having found the
solutions of the equations (\ref{MainEq}) for $\vkappa=\trr$
and/or $\vkappa=\trl$, one obtains an effective means of
studying the spectral properties of the transfer function
$M_{\Gamma_l}(z)$, referring to well known facts of operator theory
\cite{GK,Kato}.

Let $\Gamma_l$ $(l=\pm1)$ be a contour described in the hypothesis
of Lemma \ref{M1-Continuation} and let the numerical range of the operator
$\tA$ be separated from $\Gamma_l$, that is,
\begin{equation}
\label{dist-tA}
d(\Gamma_l):= \dist\big(\nu(\tA),\Gamma_l\big)>0\,.
\end{equation}
Then it is
obvious that the following quantity
\begin{equation}
\label{VarBtAGdef}
\Var_{\tA}(K,\Gamma_l) :=
\int_{\Gamma_l}|d\mu|\,\frac{\|K'(\mu)\|}{\dist\big(\mu,\nu(\tA)\big)}
\end{equation}
is finite,
\begin{eqnarray}
\nonumber
\Var_{\tA}(K,\Gamma_l) &\leq&
\mathop{\rm sup}
 \limits_{\mu\in\Gamma_l}
\frac{1+|\mu|}{\dist\big(\mu,\nu(\tA)\big)} \\
\label{VarBtAGfin}
&&\quad\times
\displaystyle\int_{\Gamma_l}|d\mu|\,(1+|\mu|)^{-1}\|K'(\mu)\|< \infty\,.
\end{eqnarray}

\begin{hypothesis}
\label{Gadmiss}
Assume all the assumptions of Section {\rm2} concerning the
operators $\tA$ and $C$, and the operator-valued function $K$.
Assume, in addition, that for both $l=+1$ and $l=-1$ there are
contours $\Gamma_l=\widetilde{\Gamma}_l\cup[\beta,+\infty)$,
described in the hypothesis of Lemma {\rm\ref{M1-Continuation}},
such that the following estimates hold true:
\begin{equation}
\label{Best}
\Var_{\tA}(K,\Gamma_l)<1\,, \qquad
\Var_{\tA}(K,\Gamma_l)\|\tA\|<
\displaystyle\frac{1}{4}\,d(\Gamma_l)\,[1-\Var_{\tA}(K,\Gamma_l)]^2\,.
\end{equation}
In the following such contours $\Gamma_l$ $(\lambda=\pm1)$ are
called admissible contours.
\end{hypothesis}

It is convenient to rewrite the equations (\ref{MainEq}),
$\vkappa=\trl$ or $\vkappa=\trr$, in the equivalent form
\begin{equation}
\label{MainEqC}
X = V^\vkappa_{\Gamma_l}(\tA+X),
\end{equation}
where $X:=Z-\tA$.
We have the following statement regarding the solvability
of the transformation equation~(\ref{MainEqC})
\begin{theorem}
\label{Solvability}
Assume Hypothesis {\rm\ref{Gadmiss}} and let $\Gamma_l$
$(l=\pm1)$ be an admissible contour.  Then for each
$\vkappa=\trr$ and $\vkappa=\trl$ the
equation~{\rm(\ref{MainEqC})} has a solution
$X_\vkappa\in\cB(\cH_A)$ satisfying the estimate
$$
   \|X_\vkappa\|\leq  r_{\rm min}(\Gamma_l),
$$
where
\begin{equation}\label{rmin}
\begin{array}{rcl}
\displaystyle r_{\rm min}(\Gamma_l)&:
=&\displaystyle \frac{1}{2}\,d(\Gamma_l)\,
[1-\Var_{\tA}(K,\Gamma_l)]\\[3mm]
&&\displaystyle -\sqrt{\frac{1}{4}\,d^2(\Gamma_l)\,[1-\Var_{\tA}(K,\Gamma_l)]^2
-d(\Gamma_l)\,\Var_{\tA}(K,\Gamma_l)\,\|\tA\|}.\quad
\end{array}
\end{equation}
Moreover, for given $\vkappa=\trr$ or $\vkappa=\trl$
this solution is a unique solution to the equation
{\rm(\ref{MainEqC})} in the open ball
$$
\bigl\{Y\in \cB(\cH_A)\, :\,
\|Y\|< r_{\rm max}(\Gamma_l)\bigr\},
$$
where
\begin{equation}\label{rmax}
r_{\rm max}(\Gamma_l)  :=  d(\Gamma_l)-\sqrt{\Var_{\tA}(K,\Gamma_l)\,
d(\Gamma_l)\,[d(\Gamma_l)+\|\tA\|]}.
\end{equation}
\end{theorem}
\begin{proof}
One can prove this theorem making use of Banach's Fixed Point
Theorem.  The proof is very similar to the proof of Theorem~3.1
in~\cite{HMM2} (also cf. \cite{MennMotMathNachr}).  Thus we omit
it.
\end{proof}
\begin{remark}
\label{Half1mVar}
Conditions (\ref{Best}) imply that the following inequalities hold
true:
$$
   r_{\rm min}(\Gamma_l) < \frac{1}{2}\,
   d(\Gamma_l)\,[1-\Var_{\tA}(K,\Gamma_l)]
    < r_{\rm max}(\Gamma_l)\,.
$$
\end{remark}
\begin{lemma}\label{Hunique}
Assume Hypothesis {\rm\ref{Gadmiss}}. Then for given $l=\pm1$ and
$\vkappa=\trr$ or $\vkappa=\trl$ the solution $X_\vkappa$
of the equation {\rm(\ref{MainEqC})} guaranteed by
Theorem {\rm\ref{Solvability}} is the same for any admissible
contour $\Gamma_l$. Moreover, this solution satisfies
the estimate
$$
\|X\|\leq r_0^{(l)}(K)
$$
where
$$
  r_0^{(l)}(K):=\inf \bigl\{r_{\rm min}(\Gamma_l)\, :\,
  \Var_{\tA}(K,\Gamma_l)<1\,,  \,\omega(K,\Gamma_l)>0\bigr\}
$$
with $r_{\rm min}(\Gamma_l)$ given by~{\rm(\ref{rmin})} and
$$
\omega(K,\Gamma_l) := d(\Gamma_l)\,[1-\Var_{\tA}(K,\Gamma_l)]^2
-4\|\tA\|\Var_{\tA}(K,\Gamma_l).
$$
\end{lemma}
\begin{proof}
This statement can be proven essentially in the same way as Theorem
3.3 in \cite{MennMotMathNachr}.
\end{proof}
Therefore, for a given holomorphy domain $\cD^{(l)}$ $(l=\pm 1)$
and fixed $\vkappa=\trr$ or $\vkappa=\trl$ the solution $X_\vkappa$
to the equation (\ref{MainEqC}) and the solution
$Z_\vkappa=\tA+X_\vkappa$ to the equation (\ref{MainEq}) do not
depend on the admissible contours $\Gamma_l$.  But when the
index $l$ changes, $X_\vkappa$ and $Z_\vkappa$ can also change.
For this reason we shall supply them in the following, when it
is necessary, with the index $l$ writing $X_\vkappa^{(l)}$ and
$Z_\vkappa^{(l)}$, $Z_\vkappa^{(l)}=\tA+X_\vkappa^{(l)}$.
Surely, the equations~(\ref{MainEq}) and~(\ref{MainEqC}) are
nonlinear equations and, outside the balls $\|X\|<r_{\rm
max}(\Gamma_l)$, they may, in principle, have other solutions,
different from the $X_\vkappa^{(l)}$ or $Z_\vkappa^{(l)}$ the
existence of which is guaranteed by Theorem~\ref{Solvability}.
In the following we only deal with the solutions
$X_\vkappa^{(l)}$ or $Z_\vkappa^{(l)}$ guaranteed by
Theorem~\ref{Solvability}.

\newsection{Factorization}\label{SecFactor}
In this section we prove a {\it factorization theorem} for the transfer
function $M_{\Gamma_l}(z)$. Note that this theorem resembles the
corresponding statements from \cite{MarkusMatsaev,VirozubMatsaev}.
It is an extension of Theorem 4.1 in \cite{HMM2}.
\begin{theorem}
\label{factorization}
Assume Hypothesis {\rm\ref{Gadmiss}} and let $\Gamma_l$ $(l=\pm1)$
be an admissible contour. Let for $\vkappa=\trr$ and $\vkappa=\trl$
the operators $X_\vkappa^{(l)}$ be the solutions of the
transformation equations {\rm(\ref{MainEqC})},
$\|X_\vkappa^{(l)}\|\leq r_0^{(l)}(K)$, and
$Z_\vkappa^{(l)}=\tA+X_\vkappa^{(l)}$. Then,
for $z\in\C\setminus \Gamma_l$, the
transfer function $M_{\Gamma_l}(z)$ admits the factorizations
\begin{equation}\label{Mfactor}
 M_{\Gamma_l}(z)=W^\trr_{\Gamma_l}(z)\,(Z_\trr^{(l)}-z)
=(Z_\trl^{(l)}-z)\,W^\trl_{\Gamma_l}(z),
\end{equation}
where $W^\vkappa_{\Gamma_l}(z)$, $\vkappa=\trr$ or
$\vkappa=\trl$, are bounded operators in $\cH_A$ which read
\begin{equation}\label{Mtild}
\begin{array}{rcl}
\displaystyle W^\trr_{\Gamma_l}(z)&=&\displaystyle I_A-\int_{\Gamma_l}
d\mu\,K'(\mu)\,(Z_\trr^{(l)}-\mu)^{-1}\\[7mm]
&&\displaystyle +z\,\int_{\Gamma_l}
d\mu\,K'(\mu)(z-\mu)^{-1}(Z_\trr^{(l)}-\mu)^{-1}
\end{array}
\end{equation}
and
\begin{equation}\label{Mtildl}
\begin{array}{rcl}
\displaystyle W^\trl_{\Gamma_l}(z)&=&\displaystyle I_A-\int_{\Gamma_l}
(Z_\trl^{(l)}-\mu)^{-1}\,K'(\mu)\,d\mu\\[7mm]
&&\displaystyle +z\,\int_{\Gamma_l}
(Z_\trl^{(l)}-\mu)^{-1}\,K'(\mu)\,d\mu
\end{array}
\end{equation}

If, in addition,
\begin{equation}
\label{distMG}
\dist\big(z,\nu(\tA)\big)\leq
d(\Gamma_l){[1-\Var_{\tA}(K,\Gamma_l)]/2},
\end{equation}
then
\begin{equation}
\label{Mtest}
\|W^\vkappa_{\Gamma_l}(z)-I_A\|<1\quad
\mbox{\it $\vkappa=\trr$ or $\vkappa=\trl$}
\end{equation}
and the operator $W^\vkappa_{\Gamma_l}(z)$ is boundedly
invertible, that is, $\|[W^\vkappa_{\Gamma_l}(z)]^{-1}\|<
\infty.$
\end{theorem}
\begin{proof}
For both $\vkappa=\trr$ and $\vkappa=\trl$ this statement
is proven in the same way. Thus, we only present the proof
for the ``right'' case $\vkappa=\trr$.

First we prove the factorization~(\ref{Mfactor}).
Note that, according to~(\ref{V1Y}) and~(\ref{MainEqC}),
\begin{equation}\label{e4.4}
\tA = Z_\trr^{(l)}-V_{\Gamma_l}(\tA+X^{(l)})
 = Z_\trr^{(l)}-\int_{\Gamma_l}
d\mu\,K'(\mu )Z_\trr^{(l)}(Z_\trr^{(l)}-\mu )^{-1}\,.
\end{equation}
Thus, in view of (\ref{Mcmpl}) and (\ref{MGamma}),
the value of $M_{\Gamma_l}(z)$ can be written as
\begin{eqnarray*}
M_{\Gamma_l}(z)&=&\tA-z+\int_{\Gamma_l}
d\mu\,K'(\mu )\frac z{z-\mu }\\
&=&Z_\trr^{(l)}-z-\int_{\Gamma_l}
d\mu\,K'(\mu)(Z_\trr^{(l)}-\mu )^{-1}\bigl(Z_\trr^{(l)}-z\bigr)\\
&&+z\int_{\Gamma_l}
d\mu\,K'(\mu)\left[ \frac {1}{z-\mu}-(Z_\trr^{(l)}-\mu )^{-1}\right]\\
&=&\bigl (Z_\trr^{(l)}-z\bigr )-\int_{\Gamma_l}
d\mu\,K'(\mu)(Z_\trr^{(l)}-\mu )^{-1}\bigl(Z_\trr^{(l)}-z\bigr)\\
&&+z\int_{\Gamma_l}
d\mu\,K'(\mu) (z-\mu)^{-1}(Z_\trr^{(l)}-\mu )^{-1} \bigl(Z_\trr^{(l)}-z\bigr).
\end{eqnarray*}
which proves the equation (\ref{Mfactor}).  The boundedness of
the operator $W^\trr_{\Gamma_l}(z)$ for $z\in\C\setminus \Gamma_l$ is
obvious.

Further, assume that the condition (\ref{distMG}) holds true.
Using the triangle inequality, it is easy to see that this
condition yields
$$
\sup\limits_{\mu\in \Gamma_l}
|z-\mu|^{-1} \leq \frac{2}{d(\Gamma_l)[1+\Var_{\tA}(K,\Gamma_l]}.
$$
For any $x\in\cH_A$, $\|x\|=1$, we have
$$
z=\lal\tA x,x\ral+\lal(z-\tA)x,x\ral
$$
and, hence,
$$
|z|\leq \|\tA\|+\Inf\limits_{\|x\|=1}|\lal(z-\tA)x,x\ral|=
\|\tA\|+\Inf\limits_{\zeta\in\nu(\tA)}|z-\zeta|=
\|\tA\|+\dist\big(z,\nu(\tA)\big).
$$
Then it follows from (\ref{distMG}) that
$$
|z|\leq\|\tA\|+\frac{1}{2}d(\Gamma_l)[1-\Var_{\tA}(K,\Gamma_l)].
$$
Meanwhile the estimate~(\ref{Arest}) implies
\begin{eqnarray}
\nonumber
\|(\tA+X^{(l)}-\mu)^{-1}\|&=&
\|[I_A+(\tA-\mu)^{-1}X_\trr^{(l)}]^{-1}(\tA-\mu)^{-1}\|\\
\label{ResEstim}
&\leq&
\frac{1}{\dist\big(\mu,\nu(\tA)\big)-\|X_\trr^{(l)}\|}\,.
\end{eqnarray}
Therefore, taking in to account the definitions of $d(\Gamma_l)$ and
$r_{\min}(\Gamma_l)$, and Remark \ref{Half1mVar} we find
\begin{equation}\label{Est00}
\left\| \int_{\Gamma_l}
d\mu\,K'(\mu)\,(Z_\trr^{(l)}-\mu)^{-1}\right\|  \leq \
\frac{2\,\Var_{\tA}(K,\Gamma_l)}{1+\Var_{\tA}(K,\Gamma_l)}\,.
\end{equation}
Using the inequality (\ref{ResEstim}) and Remark
\rref{Half1mVar} one also finds
$$
\left\|z\int_{\Gamma_l}
d\mu\,K'(\mu)\,(Z_\trr^{(l)}-\mu)^{-1}(z-\mu)^{-1}\right\|
 \leq |z|\,\frac{2\,\Var_{\tA}(K,\Gamma_l)}{1+\Var_{\tA}(K,\Gamma_l)}\,
\sup\limits_{\mu\in \Gamma_l}|z-\mu|^{-1}\,.
$$
Finally, taking into account the second assumption in (\ref{Best}),
we obtain
\begin{eqnarray*}
\|W^\trr_{\Gamma_l}(z)-I_A\| &\leq&
\frac{2\,\Var_{\tA}(K,\Gamma_l)}{1+\Var_{\tA}(K,\Gamma_l)} \\
&& +\frac{4\,\Var_{\tA}(K,\Gamma_l)\,
\left\{\|\tA\|+\frac{1}{2}\,d(\Gamma_l)
[1-\Var_{\tA}(K,\Gamma_l)]\right\}}
{d(\Gamma_l)\,[1+\Var_{\tA}(K,\Gamma_l)]^2} \\
&=& \frac{4\,\Var_{\tA}(K,\Gamma_l)\,\,[d(\Gamma_l)+\|\tA\|]}
{d(\Gamma_l)\,[1+\Var_{\tA}(K,\Gamma_l)]^2} < 1\,
\end{eqnarray*}
and, thus, if $z$ satisfies (\ref{distMG}), then
the operator $W^\trr_{\Gamma_l}(z)$ has a bounded inverse.

The proof is complete.
\end{proof}

\begin{corollary}
\label{SpHalfVic}
The spectra $\spec(Z_\trr^{(l)})$ and
$\spec(Z_\trl^{(l)})$ of the operators
$Z_\vkappa^{(l)}=\tA+X_\vkappa^{(l)}$,
$\vkappa=\trr$ or $\vkappa=\trl$, coincide and belong
to the closed $r_0^{(l)}(K)$-neighbourhood
$$
{\mathcal O}_{r^{(l)}_0(K)}(\tA):=
\bigl\{z\in\C:\,\dist\big(z,\nu(\tA)\big)
\leq r_0^{(l)}(K)\bigr\}
$$
of the numerical range $\nu(\tA)$.  Moreover, for any admissible contour
$\Gamma_l$ these spectra coincide with a subset of the spectrum
of the transfer function $M_{\Gamma_l}(\cdot)$. More precisely,
\begin{equation}
\label{e3.7new}
\spec\bigl (M_{\Gamma_l}(\cdot)\bigr)\cap {\mathcal
O}(\tA,\Gamma_l)=\spec(Z_\trr^{(l)})=\spec(Z_\trl^{(l)}),
\end{equation}
where
$$
{\mathcal O}(\tA,\Gamma_l) := \left\{z\in\C\, : \,
\dist\big(z,\nu(\tA)\big)\leq
d(\Gamma_l)\,[1-\Var_{\tA}(K,\Gamma_l)]/{2}\right\}.
$$
\end{corollary}
\begin{proof}
This statement is an immediate consequence of the factorizations
(\ref{Mfactor}) and bounded invertibility of the operators
$W^\trr_{\Gamma_l}(z)$ and $W^\trl_{\Gamma_l}(z)$ whenever
(\ref{distMG}) holds.
\end{proof}

Let us indroduce the operator
\begin{equation}\label{e4.8neu}
\Omega^{(l)} := \displaystyle\int_{\Gamma_l}
d\mu\,\mu\,\,(Z_\trl^{(l)}-\mu)^{-1}K'(\mu )\,(Z_\trr^{(l)}-\mu)^{-1},
\quad l=\pm 1,
\end{equation}
where $\Gamma_l$ stands for an admissible contour and, as
before, $Z_\vkappa^{(l)}=\tA+X_\vkappa^{(l)}$ where
$X_\vkappa^{(l)}$, $\vkappa=\trr$ or $\vkappa=\trl$, are
solutions of the transformation equations~(\ref{MainEqC}).  It
is obvious that $\Omega^{(l)}$ does not depend on the choice of
the admissible contour $\Gamma_l$.
\begin{theorem}\label{MHOmega}
The operators $\Omega^{(l)}$ $(l=\pm 1)$ possess the following properties
{\rm(}cf.
{\rm\cite{HMM1,MarkusMatsaev,MennMotMathNachr,MenShk,VirozubMatsaev}):}
\begin{eqnarray}
\label{Omest}
\|\Omega^{(l)}\|&<&1, \\[2mm]
\label{MOmega}
-\frac{1}{2\pi\ri}\int_\gamma dz\,[M_{\Gamma_l}(z)]^{-1} &=&
(I_A+\Omega^{(l)})^{-1}\,,\\[2mm]
\label{HOmega}
-\frac{1}{2\pi\ri}\int_\gamma dz\,z\,[M_{\Gamma_l}(z)]^{-1} &=&
(I_A+\Omega^{(l)})^{-1}Z_\trl^{(l)}\\
\nonumber
&=&Z_\trr^{(l)}(I_A+\Omega^{(l)})^{-1},
\end{eqnarray}
where $\gamma$ stands for an arbitrary rectifiable closed
contour encircling the spectrum
$\spec(Z_\trr^{(l)})=\spec(Z_\trr^{(l)})$ inside the set
${\mathcal O}(\tA,\Gamma_l)$ in the anticlockwise direction.
The integration along $\gamma$ is understood in the sense of the
operator norm topology.
\end{theorem}
\begin{proof}
The estimate in (\ref{Omest}) can be proved by using the
relations (\ref{Mfactor}) following the proof of the inequality
(\ref{Mtest}).  This estimate yields that the sum
$I_A+\Omega^{(l)}$ is a boundedly invertible operator in
$\cH_A$.

To prove the formula~(\ref{MOmega}) we note that by
(\ref{Mfactor}) the following factorization
holds if $z\in{\mathcal O}(\tA,\Gamma_l)\backslash \spec
(Z_\trr^{(l)})$:
\begin{equation}\label{e4.12z}
\begin{array}{rcl}
[M_{\Gamma_l}(z)]^{-1}& = &\left (Z_\trr^{(l)}-z\right)^{-1}\,
[W^\trr_{\Gamma_l}(z)]^{-1}\\[3mm]
&=&[W^\trl_{\Gamma_l}(z)]^{-1}\,
\left(Z_\trl^{(l)}-z\right)^{-1},
\end{array}
\end{equation}
where $[W^\trr_{\Gamma_l}(z)]^{-1}$ and
$[W^\trl_{\Gamma_l}(z)]^{-1}$
are holomorphic functions with values in $\cB(\cH_A)$.
By the resolvent equation and the definition (\ref{e4.8neu})
the product $\Omega^{(l)}(Z_\trr^{(l)}-z)^{-1}$
can be written as
\begin{equation}\label{e4.12zz}
\Omega^{(l)}(Z_\trr^{(l)}-z)^{-1} = F_1(z)+F_2(z),
\end{equation}
where
\begin{equation}
\label{F1int}
F_1(z) := \int_{\Gamma_l}
d\mu\,\mu\,(Z_\trl^{(l)}-\mu)^{-1}\,K'(\mu)\,(Z_\trr^{(l)}-\mu)^{-1}
\, (\mu-z) ^{-1}
\end{equation}
and
\begin{equation}\label{e4.13z}
\begin{array}{rcl}
\quad\displaystyle F_2(z)&:=&
\displaystyle\left (-\int_{\Gamma_l}d\mu
\frac {\mu}{\mu -z}\,(Z_\trl^{(l)}-\mu)^{-1}\,K'(\mu)\right )\,
(Z_\trr^{(l)}-z)^{-1}\\[3mm]
&=&
\displaystyle \left(W^\trl_{\Gamma_l}(z)-I_A\right )
(Z_\trr^{(l)}-z)^{-1}.
\end{array}
\end{equation}
Further, the formula (\ref{e4.12z})
yields
$$
(I_A+\Omega^{(l)})\,[M_{\Gamma_l}(z)]^{-1} = \
F_1(z)\,[W^\trr_{\Gamma_l}(z)]^{-1}+(Z_\trl^{(l)}-z)^{-1}.
$$
The function $F_1(z)$ is holomorphic inside the contour
$\gamma\subset{\mathcal O}(\tA,\Gamma_l)$ since the
argument $\mu$ of the integrand in (\ref{F1int}) belongs to
$\Gamma_l$ and  thereby
$$
|z-\mu|\geq [d(\Gamma_l)+\Var_{\tA}(K,\Gamma_l)]/2>0.
$$
Thus
the term $F_1(z)[W^\trr_{\Gamma_l}(z)]^{-1}$ does not contribute to
the integral
$$
-\D\frac{1}{2\pi\ri}\Int_\gamma dz (I_A+\Omega^{(l)})
[M_{\Gamma_l}(z)]^{-1}
$$
while the resolvent $(Z_\trl^{(l)}-z)^{-1}$
gives the identity $I_A$ which proves (\ref{MOmega}).

Similarly, to prove the first equality in (\ref{HOmega}) we calculate
\begin{eqnarray*}
   \lefteqn{-\D\frac{1}{2\pi\ri}\Int_\gamma dz (I_A+\Omega^{(l)})
   \,z\, [M_{\Gamma_l}(z)]^{-1} =}\\
  &=&-\D\frac{1}{2\pi\ri}\Int_\gamma dz
  \,z\, F_1(z)\,[W^\trr_{\Gamma_l}(z)]^{-1}
  -\D\frac{1}{2\pi\ri}\Int_\gamma dz \,z\,(Z_\trl^{(l)}-z)^{-1}\,.
\end{eqnarray*}
The first integral vanishes whereas the second integral equals
$Z_\trl^{(l)}$. The second equality of~(\ref{HOmega})
can be checked in the same way.
\end{proof}
\begin{corollary}\label{HlHml}
Given $l=\pm1$ the formula~{\rm(\ref{HOmega})} implies that
the operators $Z_\trr^{(l)}$ and $Z_\trl^{(l)}$ are
similar to each other:
$$
Z_\trl^{(l)} =
(I_A+\Omega^{(l)})\,Z_\trr^{(l)}\,(I_A+\Omega^{(l)})^{-1}.
$$
\end{corollary}
\begin{remark}
The formulae~(\ref{MOmega}) and~(\ref{HOmega}) allow, in
principle, to construct the operators $Z_\trr^{(l)}$ $(l=\pm 1$,
$\vkappa=\trr$ or $\vkappa=\trl)$ and, thus, to resolve the
equations~(\ref{MainEqC}) by a contour integration of the
inverse of the transfer function $M_{\Gamma_l}(z)$.
\end{remark}
\begin{theorem}
\label{MPOmega}
Let $\lambda$ be an isolated eigenvalue of
$Z_\trr^{(l)}$ and, hence, of $Z_\trl^{(l)}$ and
$M_{\Gamma_l}(\cdot)$ where $\Gamma_l$
is an admissible contour.
Denote by $P_{\trr,\lambda}^{(l)}$ and
$P_{\trl\lambda}^{(l)}$ the eigenprojections of the
operators $Z_\trr^{(l)}$ and $Z_\trl^{(l)}$, respectively, and by
$P_{M,\lambda}^{(l)}$ the residue of $M_{\Gamma_l}(z)$ at
$z=\lambda$,
\begin{equation}
\label{Plambda}
P_{\vkappa,\lambda}^{(l)} :=
-\D\frac{1}{2\pi\ri}\Int_\gamma dz\,\,(Z_\vkappa^{(l)}-z)^{-1}
\quad(\mbox{\it $\vkappa=\trr$ or $\vkappa=\trl$}),
\end{equation}
and
\begin{equation}
\label{Mpro}
P_\lambda^{(l)} := \
-\D\frac{1}{2\pi\ri}\Int_\gamma dz\,\,[M_{\Gamma_l}(z)]^{-1},
\end{equation}
where $\gamma$ stands for an arbitrary rectifiable closed
contour going around $\lambda$ in the positive direction
in a sufficiently close neighbourhood such that
$\gamma\cap\Gamma_l=\emptyset$
and no points of the spectrum of $M_{\Gamma_l}(\cdot)$, except the
eigenvalue $\lambda$, lie inside $\gamma$. Then the following
relations hold:
\begin{equation}\label{MresiduePP}
P_{M,\lambda}^{(l)} =
P_{\trr,\lambda}^{(l)}\,\,(I_A+\Omega^{(l)})^{-1}
=(I_A+\Omega^{(l)})^{-1}\,\,P_{\trl,\lambda}^{(l)}\,.
\end{equation}
\end{theorem}
\begin{proof}
Proof of this statement can be done in the same way as the proof
of the relation~(\ref{MOmega}), only the path of integration is
changed.
\end{proof}
\newsection{An example}\label{Example}
Let $\cH_A=\cH_C=L_2({\R})$ and $C=\sP^2+\lambda_C I_C$ where
$\sP=\ri\D\frac{d}{dx}$, $\lambda_C$ is some positive number,
and $I_C$ denotes the identity operator in $\cH_C$.
It is assumed that the domain $\dom(\sP)$ is the Sobolev space
$W_2^1({\R} )$ and the domain $\dom(C)$ is the Sobolev space
$W_2^2({\R})$. The spectrum of $C$ is absolutely continuous and
fills the semi-axis $[\lambda_C,+\infty)$. By the operator $A$
we understand the multiplication by a bounded complex-valued function
$a$, $Af=af$,
$f\in\cH_A$.

The operators $B$ and $D$ are defined on $\dom(B)=\dom(D)=W_2^1({\R})$
by
$$
   B=S\sP \quad\mbox{and}\quad D={\sP}Q,
$$
where $S$ and $Q$ are the multiplications by bounded
functions $s\in L_2(\R)$ and $q\in W_2^1({\R})$,
that is, $Qf=qf$ and $Sg=sg$ where $f,g{\in}L_2({\R})$.
Both $S$ and $D$ are densely defined closable operators.

Notice that $\dom(C^{1/2})=W^1_2({\R})$.  The proof of this
statement is based on the second representation theorem for
quadratic forms, see Theorem VI.2.23 in \cite{Kato}. It is
similar to the proof of Proposition~2.4 in \cite{FMM}.

Further, we assume that the functions $s$ and $q$ are
exponentially decreasing at infinity, so
that the estimates
\begin{equation}
\label{EstimB}
\bigl|s(x)\bigr|\leq c\,\exp\bigl(-\alpha_0|x|\bigr)
\quad\mbox{and}\quad
\bigl|q(x)\bigr| \leq c\,\exp\bigl(-\alpha_0|x|\bigr)\qquad (x\in{\R})
\end{equation}
hold with some $c\ge 0$ and $\alpha_0>0$.

For this example the operators $\tB$ and $\tD$ are given by
\begin{eqnarray*}
\tB &=& S\,\sP\,(\sP^2+\lambda_C I_C)^{-1/2}=
S\D\int_{\R} \frac {\mu}{(\mu ^2+\lambda_C)^{1/2}}\,
dE_{\sP}(\mu),\\
\tD&=&\overline{(\sP^2+\lambda_C I_C)^{-1/2}\,\sP\,Q}=\,\,
\D\int_{\R}\frac{\mu}{(\mu ^2+\lambda_C I_C)^{1/2}}\,
   d E_{\sP}(\mu )\, Q,
\end{eqnarray*}
where $\{E_{\sP}(\mu)\}_{\mu\in\R}$ denotes the spectral family
of the selfadjoint operator
$\sP$. Thus
\begin{eqnarray*}
\tA &=& A-\tB\tD\\
&=&\D A- S \int_{\R} \frac {\mu}{(\mu ^2+\lambda _C)^{1/2}}\,
d E_{\sP}(\mu )\,
\int_{\R} \frac {\widetilde \mu}{(\widetilde \mu ^2+\lambda _C)^{1/2}}\,
d E_{\sP}(\widetilde \mu )\, Q \\
&=&\D A- S \int_{\R} \frac {\mu^2}{\mu ^2+\lambda _C}\,
d E_{\sP}(\mu )\, Q \\
&=&A-SQ+\lambda_C S(\sP^2+\lambda_C I_C)^{-1}Q.
\end{eqnarray*}
The operator $A-SQ$ is the multiplication
by the function
$$
    \ta(x) = a(x)- s(x)q(x)
$$
while the term $S(\sP^2+\lambda_C I_C)^{-1}Q$ is a compact (even
Hilbert-Schmidt) operator in $L_2({\R})$.
Indeed, the inverse operator
$C^{-1}=(\sP^2+\lambda_C I_C)^{-1}$ is the integral operator
whose kernel reads
$$
C^{-1}(x,x') = \frac{1}{2\sqrt{\lambda_C}}\,\exp\bigl (
-\sqrt{\lambda_C}\, |x-x'|\bigr).
$$
Thus, the double integral
$\int_{\R}\int_{\R} \bigl|(S C^{-1} Q)(x,x')\bigr|^2 dx dx'$
is convergent. Obviously,
$$
\D\int_{\R}\int_{\R} \bigl |(S C^{-1} Q)(x,x')\bigr |^2\, dx\,dx' \leq
 \frac{1}{4{\lambda_C}}\,\|s\|^2_{L_2({\R})}\|q\|^2_{L_2({\R})}.
$$
Therefore, the essential spectrum of $\tA$ coincides with the range
of the function $\ta$.  In the following we assume that there
are an interval $[\alpha_1,\alpha_2]\subset(\lambda_C,+\infty)$
with $\alpha_1<\alpha_2$ and a number $\eta>0$ such that all the
numerical range $\nu(\tA)$ of $\tA$ lyes inside the domain
$$
\cO_\eta([\alpha_1,\alpha_2]):=\{z\in\C:\,
\dist(z,[\alpha_1,\alpha_2])\leq\eta\}, \quad \eta>0,
$$
of a finite real interval $[\alpha_1,\alpha_2]\subset\R$,
and, moreover, $\alpha_1-\eta>\lambda_C$.

It is easy to check that the spectral pojections $E_C(\mu)$
of the operator $C=\sP^2+\lambda_C I$ are given
by the integral operator whose kernel reads
$$
  {E}_C(\mu;x,x') = \left\{\begin{array}{cl}
0 & \mbox{ if } \mu <\lambda _C, \\[4mm]
\D\frac{1}{\sqrt{2\pi}} \int_{\lambda_C}^\mu
\D\frac{\cos[(\mu'-\lambda_C)^{1/2}(x-x')]}
{(\mu'-\lambda_C)^{1/2}}\, d\mu'
 & \mbox{ if } \mu\geq \lambda _C.
\end{array}\right.
$$
Thus, the derivative $K'(\mu)$ is also an integral operator in
$L_2({\R})$.  Its kernel $K'(\mu;x,x')$ is only nontrivial for
$\mu>\lambda_C$ and, moreover, for these $\mu$
$$
K'(\mu;x,x') =\
\D\frac{(\mu-\lambda_C)^{1/2}}{\sqrt{2\pi}\,\mu}\,
\cos[(\mu-\lambda_C)^{1/2}(x-x')]
\,\, {s(x)}\,q(x').
$$
Obviously, this kernel is degenerate for $\mu>\lambda _C$,
\begin{equation}\label{KDegenerate}
K'(\mu;x,x') = \D\frac{(\mu -\lambda _C)^{1/2}}{2\,\sqrt{2\pi}\,\mu}\,
[{s_+(\mu,x)}\,q_-(\mu,x')+s_-(\mu,x)\,q_+(\mu,x')],
\end{equation}
where $s_\pm(\mu,x)={\rm e}^{\pm\ri\,(\mu-\lambda_C)^{1/2}x}\,q(x)$
and $q_\pm(\mu,x)={\rm e}^{\pm\ri\,(\mu-\lambda_C)^{1/2}x}\,q(x)$.
{}From the assumptions (\ref{EstimB}) on $s$ and
$q$ we conclude that in the domain
$\pm \Img\sqrt{\mu-\lambda_C}<\alpha_0$, i.\,e., inside
the parabola
\begin{equation}
\label{HolDom}
\cD=\left\{\mu\in\C:\,
\Real\mu > \lambda_C-\alpha_0^2+\D\frac{1}{4\alpha_0^2}(\Img\mu)^2
\right\},
\end{equation}
the functions $s_\pm(\mu,\cdot)$ and $q_\pm(\mu,\cdot)$ are
elements of $L_2({\R})$. The function $K'(\mu)$ admits an analytic
continuation into this domain (cut along the interval
$\lambda_C-\alpha_0^2<\mu\leq\lambda_C$) as a holomorphic function
with values in $\cB(\cH_A)$ and the equation (\ref{KDegenerate})
implies that
$$
\|K'(\mu)\| \leq \D\frac{|\mu-\lambda_C|^{1/2}}{2\,\sqrt{2\pi}\,|\mu|}\,
\bigl [\|s_-(\mu,\cdot)\|\,\|q_-(\mu,\cdot)\|+
\|s_+(\mu,\cdot)\|\,\|q_+(\mu,\cdot)\|\bigr ].
$$
Obviously, for $\mu\ge\lambda_C$ we have
$\|s_\pm(\mu,\cdot)\|=\|s\|$ and
$\|q_\pm(\mu,\cdot)\|=\|q\|$.

Let us make our final assumption that
$\cD\supset\cO_\eta([\alpha_1,\alpha_2])$.  In this case one can
always choose a contour
$\Gamma=\widetilde{\Gamma}\cup[\beta,+\infty)$ where
$\beta>\alpha_2+\eta$ and the rectifiable Jordan curve
$\widetilde{\Gamma}\subset\cD\setminus\cO_\eta([\alpha_1,\alpha_2])$
results from continuous deformation of the interval
$(\lambda_C,\beta)$, the end points being
fixed.  Assume, in addition, that the functions $s$ and $q$ are
sufficiently small in the sense that the conditions (\ref{Best})
hold.  In such a case the contour $\Gamma$ is an admissible
contour (see Hypothesis \ref{Gadmiss}) and, thus, one can apply
all the statements of the Sections 3 and 4 to the corresponding
transfer function $M_\Gamma(z)$.

\begin{acknowledgements}
Support of this work by the Deutsche Forschungsgemeinschaft,
the Heisenberg--Landau Program, and the Russian Foundation for
Basic Research is gratefully acknowledged.
\end{acknowledgements}

\address{Department of Mathematics\\
University of Regensburg\\
D-93040 Regensburg\\ Germany}

\address{Department of Mathematics\\
University of Regensburg\\
D-93040 Regensburg\\ Germany}

\address{Laboratory of Theoretical Physics\\
Joint Institute for Nuclear Research\\
141980 Dubna (Moscow Region)\\
Russia}




\begin{references}{ALMSa}

\bibitem[AL]{AdL}
    {\rm Adamjan, V.\,M.,} {\rm Langer, H.}, Spectral properties
    of a class of operator-valued functions,
    {\it J.~Operator Theory} {\bf 33} (1995), 259--277.%
\bibitem[ALMSa]{AdLMSr}
    {\rm Adamyan, V.\,M., Langer, H., Mennicken, R.,} {\rm Saurer, J.,}
    Spectral components of selfadjoint block operator matrices
    with unbounded entries,
    {\it Math. Nachr.} {\bf 178} (1996), \mbox{43--80.}%
\bibitem[ALT]{AdLT}
    {\rm Adamyan, V., Langer, H.,} and {\rm Tretter, C.:}
    Existence and uniqueness of contractive solutions of some
    Riccati equations, {\it J.~Funct. Anal.} {\bf 179} (2001),
    448--473.%
\bibitem[FMM]{FMM}
  {\rm Faierman, M., Mennicken, R.,} and {\sc  M\"oller, M.}:
  The essential spectrum of a system of singular ordinary differential
  operators of mixed order. Part I: The general problem and an
  almost regular case,
  {\it Math. Nachr.} {\bf 208} (1999), 101--115.%
\bibitem[GK]{GK} {\rm Gohberg, I.\,C.,} {\rm Krein, M.\,G.},
     {\it Introduction to the theory of linear non-selfadjoint
     operators}, American Mathematical Society, Providence, 1988.%
\bibitem[HMM1]{HMM1}
     {\rm Hardt, V., Mennicken, R.,} and {\rm Motovilov, A.\,K.},
     A factorization theorem for the transfer function
     of a $2\times2$ operator matrix having
     unbounded couplings, in {\it Spectral and evolution
     problems} (Natl. Taurida Univ.
     ``V. Vernadsky", Simferopol) {\bf 10} (2000), 56--65  ({\it LANL
     e-print math.SP/9912220}).%
\bibitem[HMM2]{HMM2}
     {\rm Hardt, V., Mennicken, R.,} and {\rm Motovilov, A.\,K.},
     A factorization theorem for the transfer function
     associated with a $2\times2$ operator matrix having
     unbounded couplings, to appear in {\it J. Operator Theory}.
\bibitem[K]{Kato}
    {\rm Kato, T.}, {\it Perturbation theory for linear
    operators}, Springer-Verlag, New York,  1966.%
\bibitem[LMMT]{LMMT} Langer, H., Markus, A., Matsaev, V.,
    Tretter, C., A new concept for block oprator matrices:
    The quadratic numerical range, preprint, to appear
    in {\it Linear Algebra Appl.}%
\bibitem[LT]{LT} Langer, H., Tretter, C., Spectral decomposition
    of some nonselfadjoint block operator matrices, {\it J. Operator
    Theory} {\bf 39} (1998), 1--20.%
\bibitem[MrMt]{MarkusMatsaev}
     {\rm Markus, A.\,S.,} {\rm Matsaev, V.\,I.,}
     On the basis property for a certain part of the eigenvectors
     and associated vectors of a selfadjoint operator pencil,
     {\it Math. USSR Sb.} {\bf 61} (1988), 289--307.%
\bibitem[MM]{MennMotMathNachr}
       {\rm Mennicken, R.,} {\rm Motovilov, A.\,K.},
        Operator interpretation of resonances
        arising in spectral problems for
        ${2}\times{2}$ operator matrices,
        {\it Math. Nachr.} {\bf 201} (1999), 117--181
        ({\it LANL E-print funct-an/9708001}).%
\bibitem[MS]{MenShk}
      {\rm Mennicken, R.,} {\rm Shkalikov, A.\,A.,}
      Spectral decomposition
      of symmetric operator matrices,
      {\it Math. Nachr.} {\bf 179} (1996), 259--273.%
\bibitem[M1]{MotSPbWorkshop}
     {\rm Motovilov, A.\,K.,}
     Potentials appearing after  removal of the
     energy--dependence and scattering by them, In:
     {\it Proc. of the Intern. Workshop
     ``Mathematical aspects of the scattering
     theory and applications''},  St.~Petersburg
     University, St.~Petersburg (1991), 101--108.%
\bibitem[M2]{MotRem}
     {\rm Motovilov, A.\,K.,} Removal of the resolvent-like
     energy dependence from interactions and invariant subspaces
     of a total Hamiltonian, {\it J.~Math. Phys.} {\bf 36}  (1995),
     6647--6664 ({\it LANL E-print funct-an/9606002}).%
\bibitem[RS]{ReedSimonIII}
     {\rm Reed, M.,} {\rm Simon, B.},
     {\it Methods of modern mathematical physics,
     III: Scattering theory},  Academic Press, N.Y., 1979.%
\bibitem[VM]{VirozubMatsaev}
    {\rm Virozub, A.\,I.,} {\rm Matsaev, V.\,I.},
    The spectral properties of a certain class of selfadjoint
    operator functions,
    {\it Funct. Anal. Appl.} {\bf 8} (1974), 1--9.%
\end{references}
\end{document}